\documentstyle[11pt,amssymb]{article}

\newtheorem{prop}{Proposition}
\newtheorem{thm}{Theorem}
\newtheorem{cor}{Corollary}
\newtheorem{lemma}{Lemma}
\newtheorem{defn}{Definition}

\def\P{{\Bbb P}}

\def\Z{{\Bbb Z}}
\def\Q{{\Bbb Q}}
\def\R{{\Bbb R}}
\def\C{{\Bbb C}}
\def\O{{\cal O}}
\def\E{{\cal E}}

\let\sec\S
\let\S\relax
\def\S{{\frak S}}
\def\g{{\frak g}}
\def\p{{\frak p}}
\def\h{{\frak h}}
\def\n{{\frak n}}
\def\v{{\frak v}}

\def\k{{\frak k}}
\def\r{{\frak r}}
\def\a{{\alpha}}
\def\b{{\beta}}
\def\om{{\omega}}
\def\Om{{\Omega}}
\def\endproof{\hfill $\Box$}

\newcommand{\dis}{\displaystyle}
\newcommand{\hneq}{\hspace{-2pt}\neq\hspace{-2pt}}
\newcommand{\hin}{\hspace{-2pt}\in\hspace{-2pt}}

\newcommand{\gequ}{\geqslant}
\newcommand{\lequ}{\leqslant}
\newcommand{\ra}{\rightarrow}

\newcommand{\dbar}{\overline{\partial}}
\newcommand{\Spec}{\mbox{Spec}}

\newcommand{\End}{\mbox{End}}

\newcommand{\Ker}{\mbox{Ker}}
\newcommand{\Tr}{\mbox{Tr}}
\newcommand{\id}{\mbox{id}}
\newcommand{\Inv}{\mbox{Inv}}
\newcommand{\ov}{\overline}
\newcommand{\noin}{\noindent}
\newcommand{\wt}{\widetilde}
\newcommand{\wh}{\widehat}
\newcommand{\pr}{\prime}
\newcommand{\rk}{\mbox{rk}}

\begin{document}

\title{\bf Arithmetic Intersection Theory  on Flag Varieties}
\author{Harry Tamvakis\\ Department of Mathematics\\ University of Chicago\\
Chicago IL 60637}
\date{}
\maketitle

\begin{abstract}
 Let $F$ be the complete flag variety over $\Spec \Z$ with the
tautological filtration
$0 \subset E_1 \subset E_2 \subset\cdots \subset E_n=E$ of the
trivial bundle $E$ over $F$. The trivial hermitian metric on $E(\C)$
induces metrics on the quotient line bundles
$L_i(\C)$. Let 
$\wh{c}_1(\ov{L}_i)$ be the first Chern class of $\ov{L}_i$
in the arithmetic Chow ring $\wh{CH}(F)$ and
$\wh{x}_i=-\wh{c}_1(\ov{L}_i)$.
Let $h\hin\Z[X_1,\ldots,X_n]$ be a polynomial in the ideal
$\left<e_1,\ldots,e_n\right>$ 
generated by the elementary symmetric polynomials $e_i$.
We give an effective algorithm for computing 
the arithmetic intersection $h(\wh{x}_1,\ldots,\wh{x}_n)$ in $\wh{CH}(F)$,
as the class of a $SU(n)$-invariant differential form on $F(\C)$. 
In particular we show that all the arithmetic Chern numbers one obtains 
are rational numbers.
The results are true for partial flag varieties and 
generalize those of Maillot [Ma] for grassmannians. 
An `arithmetic Schubert calculus' is established for an `invariant
arithmetic Chow ring' which specializes to the Arakelov
Chow ring in the grassmannian case.

\end{abstract}

\section{Introduction}
\label{intro}

  Arakelov theory is a way of `completing' a variety defined over the
ring of integers of a number field by adding fibers over the archimedian
places. In this way one obtains a theory of intersection numbers using
an arithmetic degree map; these numbers are generally real valued.
The work of Arakelov on arithmetic surfaces has been generalized 
to higher dimensions by H. Gillet and C. Soul\'{e}. This provides a
link between number theory and hermitian complex geometry; the road
is via arithmetic intersection theory.

One of the difficulties with the higher dimensional theory is a lack
of examples where explicit computations are available. The arithmetic
Chow ring of projective space was studied by Gillet and Soul\'{e}
([GS2], \sec 5)  and arithmetic intersections 
on the grassmannian by Maillot [Ma].
In this article we study arithmetic intersection theory
on general flag varieties and solve two problems: (i)
finding a method to compute products in the arithmetic Chow ring,
 and (ii) formulating an `arithmetic Schubert calculus' 
analogous to the geometric case. 
The grassmannian case is easier to work with because the fiber at infinity
is a hermitian symmetric space. 
To the author's knowledge this work is the first to provide
explicit calculations when the harmonic forms are not a
subalgebra of the space of smooth forms.
The question of computing arithmetic 
intersection numbers on flag manifolds was raised by C. Soul\'{e}
in his 1995 Santa Cruz lectures [S].

We now describe our results in greater detail. The crucial case is
that of complete flags, so we discuss that for simplicity.
Let $F$ denote the complete flag variety over $\Spec\Z$,
parametrizing over any field $k$ the complete flags in a 
$k$-vector space of dimension $n$. 
Let $\ov{E}$ be the trivial vector bundle over $F$ equipped 
with a trivial hermitian metric on $E(\C)$. 
There is a tautological filtration
\[
\dis
\E:\ E_0=0 \subset E_1 \subset E_2 \subset\cdots \subset E_n=E
\]
and the metric on $E$ induces metrics on all the subbundles
$E_i$. We thus obtain a {\em hermitian filtration} $\ov{\E}$ with
quotient line bundles $L_i=E_i/E_{i-1}$, which are also given induced
metrics.
Let $\wh{CH}(F)$ be the arithmetic Chow ring of $F$ (see \sec \ref{ait} and
 [GS1], 4.2.3)
and $\wh{x}_i:=-\wh{c}_1(\ov{L}_i)$, where $\wh{c}_1(\ov{L}_i)$ is the
arithmetic first Chern class of $\ov{L}_i$ ([GS2], 2.5). 

Let $h\hin\Z[X_1,\ldots,X_n]$ be a polynomial in the ideal
$\left<e_1,\ldots,e_n\right>$ 
generated by the elementary symmetric polynomials $e_i(X_1,\ldots,X_n)$.
Our main result is a computation of
the arithmetic intersection $h(\wh{x}_1,\ldots,\wh{x}_n)$ in $\wh{CH}(F)$,
as a class corresponding to a $SU(n)$-invariant differential form on $F(\C)$. 
This enables one to reduce the computation of any
 intersection product in $\wh{CH}(F)$ to the level of smooth differential
forms; we show how to do this explicitly for products of classes
 $\wh{c}_i(\ov{E_l/E_k})$.
In particular, we obtain the following result:
Let $k_i$, $1\lequ i \lequ n$ be nonnegative integers with $\sum k_i=
\dim{F}={n \choose 2}+1$. Then the arithmetic Chern number
$
\dis
\wh{\deg}(\wh{x}_1^{k_1}\wh{x}_2^{k_2}\cdots\wh{x}_n^{k_n})
$
is a rational number. 

 Let $CH(\ov{G_d})$ be the Arakelov Chow ring ([GS1], 5.1) of the
grassmannian $G_d$ over $\Spec\Z$ parametrizing $d$-planes in
$E$, with the natural invariant K\"{a}hler metric on $G_d(\C)$.
Maillot [Ma] gave a presentation of $CH(\ov{G_d})$ and constructed an 
`arithmetic Schubert calculus' in $CH(\ov{G_d})$. There
are difficulties in extending his results to flag varieties,
mainly because the Arakelov Chow group $CH(\ov{F})$ is not a {\em subring}
of $\wh{CH}(F)$. 
To overcome this problem we define, for any partial flag variety $F$,
an `invariant arithmetic Chow ring' $\wh{CH}_{inv}(F)$. This subring of
$\wh{CH}(F)$
specializes to the Arakelov Chow ring if $F(\C)$ is a hermitian symmetric 
space.

 We extend the notion of Bott-Chern forms for 
short exact sequences to filtered bundles. These forms 
give relations in $\wh{CH}(F)$ (Theorem \ref{abc}); however
they are generally not closed forms, and thus do not represent cohomology
classes. This forces us to work on the level of differential forms
in order to calculate arithmetic intersections.
We compute the Bott-Chern forms on flag varieties $F$ by using a calculation of
the curvature matrices of homogeneous vector bundles on generalized flag
manifolds due to Griffiths and Schmid [GrS].
 One thus obtains expressions for the Bott-Chern
forms in terms of invariant forms on $F(\C)$. 

The Schubert polynomials of Lascoux and 
Sch\"{u}tzenberger provide a convenient basis 
to describe the product
structure of $\wh{CH}_{inv}(F)$. Using them we formulate an
`arithmetic Schubert calculus' for flag varieties which generalizes that
of Maillot [Ma] for grassmannians. 
However explicit general formulas are lacking, as we
 cannot do these computations using purely cohomological methods.

 This paper is organized as follows. In \sec \ref{prelim} we 
review some preliminary material on Bott-Chern forms, arithmetic
intersection theory, flag varieties and Schubert polynomials.
In \sec \ref{cbcfs} we state the main tool for computing Bott-Chern forms 
(for any characteristic class) in the case of induced metrics.
 The definition and construction 
the Bott-Chern forms associated to a hermitian filtration is the content
of \sec \ref{bcfffb}.
\sec \ref{hvb} is concerned with the explicit computation of
the curvature matrices of the tautological
vector bundles over flag varieties $F$.
In \sec \ref{afvs} we define the invariant arithmetic Chow ring
$\wh{CH}_{inv}(F)$. This subring of $\wh{CH}(F)$ is where all
the intersections of interest take place. 
 In \sec \ref{ai} we give an algorithm for calculating arithmetic intersection
numbers on the complete flag variety $F$, in particular proving that
they are all rational. In \sec \ref{asc} we describe the product structure
of $\wh{CH}_{inv}(F)$ in more detail, formulating an arithmetic Schubert
calculus. 
Some applications of our results are given in  \sec \ref{ex}. One has the
Faltings height of the image of $F$ under its pluri-Pl\"{u}cker embedding; 
this is always a rational number.  We give a table of the arithmetic Chern
numbers for $F_{1,2,3}$. Finally \sec \ref{pfvs} shows how to generalize
the previous results to partial flag varieties.

This paper will be part of the author's 1997 University of Chicago thesis.
I wish to thank my advisor William Fulton for many useful conversations and 
exchanges of ideas.

 The geometric aspects of this work generalize readily to other semisimple
groups. We plan a sequel discussing arithmetic intersection theory on
 symplectic and orthogonal flag varieties.

\section{Preliminaries}
\label{prelim}

\subsection{Bott-Chern forms}
\label{bcfs}

 The main references for this section are [BC] and [GS2].

 Consider the coordinate ring $\C[T_{ij}]$ $(1\lequ i,j\lequ n$) of the 
space $M_n(\C)$ of $n\times n$ matrices. $GL_n(\C)$ acts on matrices by 
conjugation; let $I(n)=\C[T_{ij}]^{GL_n(\C)}$ denote the 
corresponding graded ring of invariants. There is an isomorphism
$\sigma : I(n)\ra \C[X_1,X_2,\ldots,X_n]^{S_n}$ obtained by evaluating
an invariant polynomial $\phi$ on the diagonal matrix 
diag$(X_1,\ldots,X_n)$. We will often identify $\phi$ with the the
symmetric polynomial $\sigma(\phi)$. 
We let $I(n,\Q)=\sigma^{-1}(\Q[X_1,X_2,\ldots,X_n]^{S_n})$. For $A$
 an abelian group, $A_{\Q}$ denotes $A\otimes_{\Z}\Q$.

 Let $X$ be a complex manifold, and denote by $A^{p,q}(X)$ the space of 
differential forms of type $(p,q)$ on $X$. Let 
 $A(X)=\bigoplus_p A^{p,p}(X)$ and $\wt{A}(X)$ be the
quotient of $A(X)$ by $\mbox{Im} \partial + \mbox{Im} \dbar$. 
If $\omega$ is a closed form in $A(X)$ the cup product $\wedge\omega:
\wt{A}(X)\ra\wt{A}(X)$ and the operator
$dd^c:\wt{A}(X)\ra A(X)$ are well defined.

   Let $E$ be a rank $n$ holomorphic vector
bundle over $X$, equipped with a hermitian metric $h$. The pair 
$\ov{E}=(E,h)$ is called a {\em hermitian vector bundle}. 
A direct sum $\ov{E}_1\bigoplus\ov{E}_2$ of 
hermitian vector bundles will always
mean the orthogonal direct sum
$(E_1\bigoplus E_2,h_1\oplus h_2)$.
Let $D$ be the hermitian holomorphic connection of $\ov{E}$,
with curvature $K=D^2\hin A^{1,1}(X,\End(E))$.

  If $\phi\hin I(n)$ is any invariant polynomial,
there is an associated differential form 
$\phi(\ov{E}):=\phi(\frac{i}{2\pi}K)$,
 defined locally by identifying $\End(E)$ with $M_n(\C)$.
 These differential forms are $d$ and $d^c$ closed, have de Rham cohomology 
class independent of the metric $h$, and are functorial under pull back by
 holomorphic maps (cf. [BC]). In particular one
 obtains the {\em power sum forms} $p_k(\ov{E})$ with 
$\dis p_k=\sum_i X_i^k$  and the
{\em Chern forms} $c_k(\ov{E})$ with $c_k=e_k$
the $k$-th elementary symmetric polynomial.

 Let $\E:\ 0\ra S \ra E\ra Q \ra 0$ be an exact 
sequence of holomorphic vector bundles on $X$. Choose arbitrary hermitian
metrics $h_S,h_E,h_Q$ on $S,E,Q$ respectively.
Let
\begin{equation}
\label{ses}
\ov{\E}=(\E,h_S,h_E,h_Q):\
0\ra \ov{S} \ra \ov{E}\ra \ov{Q}\ra 0.
\end{equation}
We say that $\ov{\E}$ is {\em split}
when $(E,h_E)=(S\bigoplus Q,h_S\oplus h_Q)$ and 
$\E$ is the obvious exact sequence. 

Let $\phi\hin I(n)$ be any invariant polynomial.
Then there is a unique way to attach to every exact
sequence $\ov{\E}$ a form $\wt{\phi}(\ov{\E})$ in 
$\wt{A}(X)$, called the Bott-Chern form of $\ov{\E}$, in such a way that:

\noin
(i) $dd^c\wt{\phi}(\ov{\E})=\phi(\ov{S}\bigoplus \ov{Q})
-\phi(\ov{E})$,

\noin
(ii) For every map $f:X\ra Y$ of complex manifolds,
        $\wt{\phi}(f^*(\ov{\E}))=f^*\wt{\phi}(\ov{\E})$,

\noin
(iii) If $\ov{\E}$ is split, then $\wt{\phi}(\ov{\E})=0$.

 For $\phi$, $\psi \hin I(n)$ one has the following
 useful relations in $\wt{A}(X)$:
\begin{equation}
\label{define}
\wt{\phi + \psi}=\wt{\phi}+\wt{\psi}, \ \ \ \ \ \ 
\wt{\phi\psi}=\wt{\phi}\cdot\psi(\ov{S}\oplus \ov{Q})+
\phi(\ov{E})\cdot\wt{\psi}.
\end{equation}

\subsection{Arithmetic intersection theory}
\label{ait}

 We recall here the generalization of Arakelov
theory to higher dimensions due to H. Gillet and C. Soul\'{e}. 
For more details see [GS1], [GS2], [SABK].

Let $X$ be an
{\em arithmetic scheme}, by which we mean a regular scheme, 
projective and flat over $\mbox{Spec}\Z$. 
For $p\gequ 0$, we denote the Chow group of 
codimension $p$ cycles on
$X$ modulo rational equivalence by
$CH^p(X)$ and let $CH(X)=\bigoplus_p CH^p(X)$. $\wh{CH}^p(X)$
will denote the $p$-th arithmetic Chow group of $X$. Recall that
an element of $\wh{CH}^p(X)$ is represented by an arithmetic cycle
$(Z,g_Z)$; here $g_Z$ is a Green current for the codimension $p$ cycle
$Z(\C)$. Let $\wh{CH}(X)=\bigoplus_p \wh{CH}^p(X)$.

The involution of $X(\C)$ induced by complex conjugation is denoted 
by $F_{\infty}$. Let $A^{p,p}(X_{\R})$
be the subspace of $A^{p,p}(X(\C))$ generated
by real forms $\eta$ such that $F^*_{\infty}\eta=(-1)^p\eta$;
denote by $\wt{A}^{p,p}(X_{\R})$
the image of $A^{p,p}(X_{\R})$ in $\wt{A}^{p,p}(X(\C))$. 
Let $A(X_{\R})=\bigoplus_p A^{p,p}(X_{\R})$ and 
$\wt{A}(X_{\R})=\bigoplus_p \wt{A}^{p,p}(X_{\R})$.
We have the following canonical morphisms of abelian groups:
\[
\dis
 \zeta :\wh{CH}^p(X)  \longrightarrow  CH^p(X), \ \ \ 
 {[(Z,g_Z)]} \longmapsto  {[Z]},
\]
\[
\dis
 \omega : \wh{CH}^p(X)  \longrightarrow A^{p,p}(X_{\R}), \ \ \ 
 {[(Z,g_Z)]}  \longmapsto dd^cg_Z+\delta_{Z(\C)},
\]
\[
\dis
 a : \wt{A}^{p-1,p-1}(X_{\R}) \longrightarrow \wh{CH}^p(X), \ \ \
 \eta  \longmapsto  {[(0,\eta)]}.
\]

For convenience of notation, when we refer to a real differential form
$\eta\hin A(X_{\R})$ as an element of $\wh{CH}(X)$, we shall always mean
$a([\eta])$, where $[\eta]$ is the class of $\eta$ in $\wt{A}(X_{\R})$.
There is an exact sequence
\begin{equation}
\label{ex1}
CH^{p,p-1}(X) \longrightarrow \wt{A}^{p-1,p-1}(X_{\R})
\stackrel{a}\longrightarrow \wh{CH}^p(X) 
\stackrel{\zeta}\longrightarrow CH^p(X)\longrightarrow 0
\end{equation}
Here the group $CH^{p,p-1}(X)$ is the 
$E_2^{p,1-p}$ term of a spectral sequence used by Quillen to 
calculate the higher algebraic $K$-theory of $X$ (cf. [G]).

One can define a pairing $\wh{CH}^p(X)\otimes\wh{CH}^q(X)
\ra \wh{CH}^{p+q}(X)_{\Q}$ which turns $\wh{CH}(X)_{\Q}$
into a commutative graded unitary $\Q$-algebra. The maps $\zeta$, $\omega$
are $\Q$-algebra homomorphisms. If $X$ is smooth over $\Z$
one does not have to tensor with $\Q$. The functor $\wh{CH}^p(X)$ is 
contravariant in $X$, and covariant for proper maps which are smooth on
the generic fiber.
 We also note the useful
identity $a(x)y=a(x\omega(y))$ for $x,y\hin\wh{CH}(X)$.

 Choose a K\"{a}hler form $\om_0$ on $X(\C)$ such that 
$F^*_{\infty}\omega_0=-\om_0$ and let 
$\cal{H}^{p,p}(X_{\R})$ be the space of harmonic (with 
respect to $\omega_0$) $(p,p)$ 
forms on $X(\C)$ invariant under $F_{\infty}$. The $p$-th {\em Arakelov Chow
group} of $\ov{X}=(X,\om_0)$ is defined by $CH^p(\ov{X}):=
\omega^{-1}(\cal{H}^{p,p}(X_{\R}))$. The group $CH(\ov{X})_{\Q}=\bigoplus_p
CH^p(\ov{X})_{\Q}$ is generally not a subring of $\wh{CH}(X)_{\Q}$,
unless the harmonic forms $\cal{H}^*(X_{\R})$ are a subring of
$A(X_{\R})$. This is true if $(X(\C),\om_0)$ is a hermitian symmetric
space, such as a complex grassmannian,
but fails for more general flag varieties.

Let $f:X\ra\Spec\Z$ be the projection.
  If $X$ has relative dimension $d$ over $\Z$, then we have an
arithmetic degree map $\wh{\deg}:\wh{CH}^{d+1}(X)\ra\R$, obtained by 
composing the push-forward $f_*:\wh{CH}^{d+1}(X)\ra\wh{CH}^1(\Spec\Z)$
with the isomorphism $\wh{CH}^1(\Spec\Z)\stackrel{\sim}\ra\R$. 
The latter maps the class of $(0,2\lambda)$ to the real number $\lambda$.

 A {\em hermitian vector bundle} $\ov{E}=(E,h)$ on an arithmetic scheme $X$ is
an algebraic vector bundle $E$ on $X$ such that the induced holomorphic
vector bundle $E(\C)$ on $X(\C)$ has a hermitian metric $h$ with
 $F_{\infty}^*(h)=h$. There are characteristic classes 
$\wh{\phi}(\ov{E})\hin \wh{CH}(X)_{\Q}$ for any $\phi\hin I(n,\Q)$, 
where $n=\rk E$. For example, we have {\em arithmetic Chern classes} 
$\wh{c}_k(\ov{E}) \hin \wh{CH}^k(X)$. 
 For the basic properties of these classes, see [GS2], Theorem 4.1.

\subsection{Flag varieties and Schubert polynomials}
\label{classgps}

 Let $k$ be a field, $E$ an $n$-dimensional vector space 
over $k$.
Let 
\[
\dis 
\r=(0<r_1<r_2<\ldots<r_m=n)
\]
be an increasing $m$-tuple of
natural numbers. A {\em flag of type $r$} is a flag
\begin{equation}
\label{fil}
\E:\ E_0=0 \subset E_1 \subset E_2 \subset\cdots \subset E_m=E
\end{equation}
with $\mbox{rank}{E_i}=r_i$, $1\lequ i\lequ m$. 
Let $F(\r)$ denote the arithmetic scheme parametrizing flags $\E$ 
of type $\r$ over any field $k$. (\ref{fil}) will also denote 
the tautological flag of vector bundles over $F(\r)$, and we call the 
resulting filtration of the bundle $E$ a {\em filtration of type $\r$}.

 The above {\em arithmetic flag variety} is smooth over
$\Spec\Z$.  There is an isomorphism
$F(\r)(\C)\simeq SL(n,\C)/P$, where $P$ is the parabolic subgroup 
of $SL(n,\C)$ stabilizing a fixed flag.
 In the extreme case $m=2$ (resp. $m=n$) $F(\r)$
is the grassmannian $G_d$ parametrizing $d$-planes in $E$ 
(resp. the complete flag variety $F$). Although the results of this paper 
are true for any partial flag variety $F(\r)$, for simplicity we will
work with the complete flag variety $F$, leaving the discussion of
the general case to \sec \ref{pfvs}. 
 The notation for these varieties
and the dimension $n$ will be fixed throughout this paper.

  We now  recall the standard presentation of the Chow group $CH(F)$.
Define the quotient line bundles $L_i=E_i/E_{i-1}$.
Consider the polynomial ring $P_n=\Z[X_1,\ldots,X_n]$ and the ideal
$I_n$ generated by the elementary symmetric functions $e_i(X_1,\ldots,X_n)$.
Then $CH(F)\simeq P_n/I_n$,
where the inverse of this isomorphism sends $[X_i]$ to $-c_1(L_i)$.
 The ring $H_n=P_n/I_n$
has a free $\Z$-basis consisting of classes of monomials 
$X_1^{k_1}X_2^{k_2}\cdots X_n^{k_n}$, where the exponents
$k_i$ satisfy $k_i\lequ n-i$. 

The Schubert polynomials of Lascoux and Sch\"{u}tzenberger [LS]
are a natural $\Z$-basis of $H_n$, corresponding to
the classes of Schubert varieties in $CH(F)$. Our main reference
for Schubert polynomials will be Macdonald's notes [M].

Let $S_{\infty}=\cup_nS_n$ and $P_{\infty}=\Z[X_1,X_2,\ldots]$. 
For each $w\hin S_{\infty}$, $l(w)$ denotes the 
{\em length} of $w$ and $\partial_w:P_{\infty}\ra P_{\infty}$
 the corresponding {\em divided difference operator}\ ([M] Chp. 2). 
 If $w_0$ is the longest element of $S_n$ and $w\hin S_n$ is arbitrary,
the {\em Schubert polynomial} $\S_w$ is given by 
\[
\dis
\S_w=\partial_{w^{-1}w_0}(X_1^{n-1}X_2^{n-2}\cdots X_{n-1}).
\]
This definition is compatible with the natural inclusion $S_n\subset
S_{n+1}$ (with $w(n+1)=n+1$). It follows that $\S_w$ is well defined for 
any $w\hin S_{\infty}$.

We let $\Lambda_n=P_n^{S_n}$ be the ring
of symmetric polynomials. 
The set $\{\S_w\ | \ w\hin S_n \}$ is both a free $\Lambda_n$-basis of
$P_n$ and a free $\Z$-basis of $H_n$. Let $S^{(n)}$ denote the set of
permutations $w\hin S_{\infty}$ such that $w(n+1)<w(n+2)<\cdots$
Then $\{\S_w\ | \ w\hin S^{(n)} \}$ is a free $\Z$-basis of $P_n$\
([M], (4.13)).
 
Define a $\Lambda_n$-valued scalar product on $P_n$ by
$\left< f ,g \right>=\partial_{w_0}(fg)$, for $f,g\hin P_n$.
If $\{\S^w\}_{w\in S_n}$ is the $\Lambda_n$-basis of $P_n$ dual
to the basis $\{\S_w\}$ relative to this product, then
$\S^w(X)=w_0\S_{ww_0}(-X)$. \ ([M], (5.12)). For any $h\hin I_n$ we
have a decomposition $\dis h=\sum_{w\in S_n}\left<h,\S^w\right>\S_w$,
where each $\left<h,\S^w\right>$ is in $\Lambda_n\cap I_n$.

\section{Calculating Bott-Chern forms}
\label{cbcfs}

 Consider the short exact sequence $\ov{\E}$ in (\ref{ses}) and assume that
the metrics on $\ov{S}$ and $\ov{Q}$ are induced from the
metric on $E$. Let $r$, $n$ be the ranks of the bundles $S$ and $E$.

For $\phi\hin I(n)$ homogeneous of degree $k$
we let $\phi^{\pr}$ be a $k$-multilinear invariant form on $M_n(\C)$ such that
$\phi(A)=\phi^{\pr}(A,A,\ldots,A)$. 
Such forms are most easily constructed for the power sums $p_k$,
 by defining
\[
p_k^{\pr}(A_1,A_2,\ldots,A_k)=\Tr (A_1A_2\cdots A_k).
\]
If $\lambda=(\lambda_1,\lambda_2,\ldots,\lambda_m)$ is a partition
of $k$, define $\dis p_{\lambda}:=\prod_{i=1}^mp_{\lambda_i}$. For 
$p_{\lambda}$ we can take $p_{\lambda}^{\pr}=\prod p_{\lambda_i}^{\pr}$.
Since the $p_{\lambda}$'s are an additive $\Q$-basis for the ring of
symmetric polynomials, we can use the above constructions to find 
multilinear forms $\phi^{\pr}$ for any $\phi\hin I(n)$.
For any two matrices $A,B \hin M_n(\C)$ let
\[
\dis
\phi^{\pr}(A;B):=\sum_{i=1}^k\phi^{\pr}(A,A,\ldots,A,B_{(i)},A,\ldots,A),
\]
where the index $i$ means that $B$ is in the $i$-th position.

Consider a local orthonormal frame $s$ for $E$ such that the first $r$ elements
generate $S$, and
let $K(\ov{S})$, $K(\ov{E})$ and $K(\ov{Q})$ 
be the curvature matrices of $\ov{S}$, $\ov{Q}$
and $\ov{E}$ with respect to $s$.
Let $K_S=\frac{i}{2\pi}K(\ov{S})$,
$K_E=\frac{i}{2\pi}K(\ov{E})$ and
$K_Q=\frac{i}{2\pi}K(\ov{Q})$. Write
\[
\dis
K_E=
\left(
\begin{array}{c|c}
K_{11} & K_{12} \\ \hline
K_{21} & K_{22} 
\end{array} \right)
\]
where $K_{11}$ is an $r\times r$ submatrix,
and consider the matrices
\[
\dis
K_0=
\left(
\begin{array}{c|c}
K_S & 0 \\ \hline
K_{21} & K_Q 
\end{array} \right)
\ \mbox{ and } \
J_r=
\left(
\begin{array}{c|c}
Id_r & 0 \\ \hline
0 & 0 
\end{array} \right).
\]
   
Let $u$ be a variable and define $K(u)=uK_E+(1-u)K_0$.
We can then state the main computational 
\begin{prop}
\label{calc} For $\phi\hin I(n)$, we have
\[
\dis
\wt{\phi}(\ov{\E})=\int_0^1\frac{\phi^{\pr}(K(u); J_r)-
\phi^{\pr}(K(0); J_r)}{u}\,du.
\]
\end{prop}
Proposition \ref{calc} is essentially a consequence of the work of
Bott and Chern [BC], although we have not been able to find this general 
statement in the literature. For history and a complete proof, see [T].

\medskip

What will prove most useful to us in the sequel is
\begin{cor}
\label{ratcor}
For any $\phi\hin I(n,\Q)$ the Bott-Chern
form $\wt{\phi}(\ov{\E})$ is a polynomial in the entries of the matrices
$K_E$, $K_S$ and $K_Q$ with {\em rational} coefficients.
\end{cor}
{\bf Proof.} By the equations (\ref{define}) it suffices to prove this 
for $\phi=p_k$ a power sum.
Using the bilinear form $p_k^{\pr}$ described previously in Proposition
\ref{calc} gives
\[
\dis
\wt{p}_k(\ov{\E})=k\int_0^1\frac{1}{u}\Tr((K(u)^{k-1}-K(0)^{k-1})J_r)\,du,
\]
so the result is clear.
\endproof

\medskip

 Define the {\em harmonic numbers}
$\dis \cal{H}_k=\sum_{i=1}^k\frac{1}{i}$, $\cal{H}_0=0$. 
We will need the following useful calculations,
which one can deduce from the definition and from Proposition \ref{calc}:

\noin
(a) $\wt{c_1^k}(\ov{\E})=0$ for all $k$ 
and $\wt{c}_p(\ov{\E})=0$ for all $p > \rk E$.

\noin
(b) $\wt{c}_2(\ov{\E})=c_1(\ov{S})-\Tr K_{11}$ (see [D], 10.1 and [T]).

\noin
(c) If $E$ is flat, then 
$\dis \wt{c}_k(\ov{\E})=
\cal{H}_{k-1}\sum_{i=0}^{k-1}ic_i(\ov{S})c_{k-1-i}(\ov{Q})$
([Ma], Th. 3.4.1).

\section{Bott-Chern forms for filtered bundles}
\label{bcfffb}

In this section we will extend the definition of Bott-Chern forms 
for an exact sequence of bundles to the case of a filtered bundle.

Let $X$ and $E$  be  as in \sec \ref{bcfs}, and assume that $E$ has
a filtration of type $\r$
\begin{equation}
\label{fil2}
\E:\ E_0=0 \subset E_1 \subset E_2 \subset\cdots \subset E_m=E
\end{equation}
by complex subbundles $E_i$, with $\r$ as in \sec \ref{classgps}.
Let $Q_i=E_i/E_{i-1}$, $1\lequ i\lequ m$ be the quotient bundles.
A {\em hermitian filtration $\ov{\E}$ of type $\r$}
is a filtration (\ref{fil2}) together with a choice of hermitian metrics
 on $E$ and on each quotient bundle $Q_i$. Note that we do not assume that
metrics
have been chosen on the subbundles $E_1,\ldots, E_m$ or that the metrics
on the quotients are induced from $E$ in any way. We say that $\ov{\E}$ is
{\em split} if, when $E_i$ is given the induced metric from $E$ for each $i$,
the sequence $\ov{\E_i} : 0\ra \ov{E}_{i-1} \ra \ov{E}_i
\ra \ov{Q}_i \ra 0$ is split. In this case of course
$\dis \ov{E}=\bigoplus_i\ov{Q}_i$.

\begin{thm} \label{bcf}
 Let $\phi\hin I(n)$ be an invariant polynomial.
 There is a unique way to attach to every hermitian
filtration of type $\r$ a form $\wt{\phi}(\ov{\E})$ in 
$\wt{A}(X)$ in such a way that:

\noin
{\em (i)} $\dis dd^c\wt{\phi}(\ov{\E})=
\phi(\bigoplus_{i=1}^m \ov{Q}_i)-\phi(\ov{E})$,

\noin
{\em (ii)} For every map $f:X\ra Y$ of complex manifolds,
        $\wt{\phi}(f^*(\ov{\E}))=f^*\wt{\phi}(\ov{\E})$,

\noin
{\em (iii)} If $\ov{\E}$ is {\em split}, then $\wt{\phi}(\ov{\E})=0$.

If $m=2$, i.e. the filtration $\E$ has length 2, then 
$\wt{\phi}(\ov{\E})$ coincides with the Bott-Chern class 
 $\wt{\phi}(0 \ra \ov{Q}_1 \ra \ov{E} \ra \ov{Q}_2 \ra 0)$
 defined in {\em \sec \ref{bcfs}}.

\end{thm}
{\bf Proof.} The essential ideas are contained in [GS2], Th. 1.2.2 and
sections 7.1.1, 7.1.2., so we will only sketch the argument.

We first show that such forms exist. 
Given any hermitian filtration
$\E$, equip each subbundle $E_i$ with the induced metric from $\ov{E}$
and consider the exact sequence 
\[
\dis
\ov{\E_i} : 0\ra \ov{E}_{i-1} \ra \ov{E}_i
\ra \ov{Q}_i \ra 0.
\]

If $\wt{\phi}(\ov{\E})$ and $\wt{\psi}(\ov{\E})$ have already been 
defined then the equations
\[
\dis
\wt{\phi + \psi}(\ov{\E})=\wt{\phi}(\ov{\E})+\wt{\psi}(\ov{\E})
\]
\[
\wt{\phi\psi}(\ov{\E})=\wt{\phi}(\ov{\E})\psi(\bigoplus_{i=1}^m \ov{Q}_i)+
\phi(\ov{E})\wt{\psi}(\ov{\E})
\]
can be used to define $\wt{\phi + \psi}$ and $\wt{\phi\psi}$ (see [GS2],
Prop. 1.3.1 for the case $m=2$). Therefore
it suffices to construct the Bott-Chern classes $\wt{p_k}$ for the
power sums $p_k$. 
For this we simply let
\begin{equation}
\label{*}
\wt{p_k}(\ov{\E}):=\sum_{i=1}^m\wt{p_k}(\ov{\E_i}).
\end{equation}
Since the $\wt{p_k}(\ov{\E_i})$ are functorial and 
additive on orthogonal direct sums, it is clear
that (\ref{*}) satisfies (i)-(iii). The construction
for $m=2$ gives the classes of \sec \ref{bcfs}.

 We will use a separate construction of
the total Chern forms $\wt{c}(\ov{\E})$: 
For each $i$ with $1\lequ i\lequ m-1$, let 
$\ov{\cal{Q}}_i$ be the sequence $\dis 0 \ra 0 \ra \bigoplus_{j=i+1}^m
\ov{Q}_j \ra \bigoplus_{j=i+1}^m \ov{Q}_j \ra 0$, and let 
$\ov{\E_i^+}=\ov{\E_i}\bigoplus \ov{\cal{Q}}_i$. Let 
$\ov{\E_m^+}=\ov{\E_m}$. 
To each exact sequence $\ov{\E_i^+}$ we associate 
a Bott-Chern form $\wt{c}(\ov{\E_i^+})$. It follows from [GS2], Prop. 1.3.2 
that
\[
\dis
\wt{c}(\ov{\E_i^+})=\wt{c}(\ov{\E_i}\bigoplus \ov{\cal{Q}}_i)
=\wt{c}(\ov{\E_i})c(\bigoplus_{j=i+1}^m\ov{Q_j})
=\wt{c}(\ov{\E_i})\wedge\bigwedge_{j=i+1}^mc(\ov{Q_j}).
\]
It is easy to see that $\dis \wt{c}(\ov{\E}):=
\sum_{i=1}^m \wt{c}(\ov{\E_i^+})$ satisfies (i)-(iii).

 To prove that the form  $\wt{\phi}(\ov{\E})$ is unique, one constructs
a deformation of the filtration $\ov{\E}$ to the split filtration,
as in [GS2], \sec 7.1.2. 
Let $\ov{\O(1)}$ be the canonical line bundle on $\P^1=\P^1(\C)$ 
with its natural Fubini-Study metric and let $\sigma$ be a 
section of $\O(1)$ vanishing only at $\infty$. Let $p_1$, $p_2$
be the projections from $X\times\P^1$ to $X$, $\P^1$ respectively. We 
denote by $E$, $E_i$, $Q_i$ and $\O(1)$ the bundles $p_1^*E$, $p_1^*E_i$,
$p_1^*Q_i$, and $p_2^*\O(1)$ on $X\times\P^1$. For a bundle $F$ on 
$X\times\P^1$ we let $F(k):=F\otimes\O(1)^k$. 

 For each $i\leqslant m-1$, we map $E_i(m-1-i)$ to $E_{i+1}(m-1-i)$ by the
inclusion of $E_{i}\hookrightarrow E_{i+1}$ and to $E_i(m-i)$ by 
$\id_{E_i}\otimes\sigma$. For $1\lequ j\lequ m$ let
\[
\dis
\wt{E}_k:=\left( \bigoplus_{i=1}^k E_i(m-i)\right)\Bigl/
\left( \bigoplus_{i=1}^{k-1} E_i(m-1-i)\right).
\]
 
 Setting $\wt{E}:=\wt{E}_m$ we get a filtration of type 
$\r$ over $X\times \P^1$ :
\[
\wt{\E} :\ 0 \subset \wt{E}_1 \subset \wt{E}_2 \subset
\cdots \subset \wt{E}_m=\wt{E}.
\]
The quotients of this filtration are
$\wt{Q}_i=\wt{E}_i/\wt{E}_{i-1}=Q_i(m-i)$,
 $1\lequ i\lequ m$.

For $z\hin\P^1$, denote by $i_z:X\ra X\times \P^1$ the map given by
$i_z(x)=(x,z)$. When $z\hneq\infty$, $i_z^*\wt{E}\simeq E$, while 
$\dis i_{\infty}^*\wt{E}\simeq
\bigoplus_{i=1}^m Q_i$.
Using a partition of unity we can choose
hermitian metrics $\wt{h}_i$ on $\wt{Q}_i$ and 
$\wt{h}$ on $\wt{E}$ such that the isomorphisms $i_z^*\wt{Q}_i\simeq Q_i$,
 $i_0^*\wt{E}\simeq E$ and $\dis i_{\infty}^*\wt{E}\simeq
\bigoplus_{i=1}^m Q_i$ all become isometries.
 We also let $(\wt{\E},\wt{h})$ denote the
hermitian filtration of type $\r$ defined by these data.
Then one shows (as in loc. cit.) that $\wt{\phi}(\ov{\E})$ is
uniquely determined in $\wt{A}(X)$ by the formula
\[
\dis
\wt{\phi}(\ov{\E})=-\int_{\P^1}\phi(\wt{E},\wt{h})\log |z|^2. 
\]
\endproof

\medskip

\noindent
{\bf Remark.} Gillet and Soul\'{e} used $\wt{ch}(\ov{\E})$
to give an explicit description of the Beilinson regulator map
on $K_1(X)$, where $X$ is an arithmetic scheme ([GS2], 7.1).

\medskip

 It is easy to prove that analogues of the properties of Bott-Chern forms
for short exact sequences ([GS2], \sec 1.3) are true for the above
generalization to filtered bundles. In particular the formulas 
(\ref{define}) take the form:

\begin{equation}
\label{sumprod}
\wt{\phi + \psi}=\wt{\phi}+\wt{\psi}, \ \ \ \ \ \ 
\wt{\phi\psi}=\wt{\phi}\cdot\psi(\bigoplus_{i=1}^m \ov{Q}_i)+
\phi(\ov{E})\cdot\wt{\psi}.
\end{equation}
for any $\phi$, $\psi \hin I(n)$.
Using Theorem \ref{bcf} and the same argument
as in the proof of Theorem 4.8(ii) in [GS2], we obtain
\begin{thm} \label{abc}
 Let
\[
\dis
\ov{\E} :\ 0\subset \ov{E}_1\subset \ov{E}_2 \subset\cdots\subset 
\ov{E}_m=\ov{E}
\]
be a hermitian filtration on an arithmetic scheme $X$, with
quotient bundles $Q_i$, and let $\phi\hin I(n,\Q)$. Then
\[
\dis
\wh{\phi}(\bigoplus_{i=1}^m\ov{Q}_i)-
\wh{\phi}(\ov{E})=a(\wt{\phi}(\ov{\E})). 
\]
\end{thm}

Assume that the subbundles $E_i$ are given metrics induced from $E$ and 
the quotient bundles $Q_i$ are given the metrics induced
from the exact sequences $\ov{\E_i}$. Define matrices 
$K_{E_i}=\frac{i}{2\pi}K(\ov{E}_i)$ and $K_{Q_i}=\frac{i}{2\pi}
K(\ov{Q}_i)$ as in \sec \ref{cbcfs}. Then the constructions in Theorem 
\ref{bcf} and Corollary \ref{ratcor} immediately imply

\begin{cor}
\label{ratbcf}
For any $\phi\hin I(n,\Q)$ the Bott-Chern
form $\wt{\phi}(\ov{\E})$ is a polynomial in the entries of the matrices
$K_{E_i}$ and $K_{Q_i}$, $1\lequ i \lequ m$, with {\em rational} coefficients.
\end{cor}

\section{Curvature of homogeneous vector bundles}
\label{hvb}

  Let $G=SL(n,\C)$, $K=SU(n)$ and $P$ be a parabolic subgroup of $G$, with
Lie algebras $\g$, $\k$ and $\p$ respectively. Complex conjugation of $\g$
with respect to $\k$ is given by the map $\tau$ with $\tau(A)=-\ov{A}^t$.
We let
$\v=\p\cap\tau(\p)$ and $\n$ be the unique maximal nilpotent ideal
of $\p$, so that $\g=\v\oplus\n\oplus\tau(\n)$. 

Let $\h=\{\mbox{diag}(z_1,\ldots,z_n)\ |\ \sum z_i=0\}$ be the Cartan
subalgebra of diagonal matrices in $\g$. The set of roots
$\Delta=\{z_i-z_j\ |\ 1\lequ i\neq j \lequ n \}$ is a subset of $\h^*$.
We denote the root $z_i-z_j$ by the pair $ij$, and fix a system of
positive roots $\Delta_+=\{ij\ |\ i>j\}$. The adjoint representation of $\h$ 
on $\g$ determines a decomposition 
$\dis \g=\h\oplus\sum_{\a\in\Delta}\g^{\a}$. Here the root space 
$\g^{\a}=\C e_{\a}$, where $e_{\a}=e_{ij}=E_{ij}$ 
is the matrix with 1 at the $ij$-th entry and zeroes elsewhere.
Set $\ov{e}_{ij}=\tau(e_{ij})=-E_{ji}$.

 Let $V=K\cap P$ and consider the complex manifold
$X=G/P=K/V$. Let $p:K\ra X$ be the quotient map, and
let $\Psi\subset\Delta_+$ be such that $\dis \n=\sum_{\a\in\Psi}\g^{-\a}$.
For $\a,\b\hin\Psi$, the equations
\begin{eqnarray*}
\omega^{\a}(e_{\b})=\delta_{\a\b}, & \omega^{\a}(\ov{e}_{\b})=0,
& \omega^{\a}(\v)=0 \\
\ov{\omega}^{\a}(e_{\b})=0, & \ov{\omega}^{\a}(\ov{e}_{\b})=\delta_{\a\b},
& \ov{\omega}^{\a}(\v)=0
\end{eqnarray*}
define elements of the dual space $\g^*$, which we shall regard as left
invariant complex one-forms on $K$. A given differential form $\eta$ on
$X$ pulls back to 
\begin{equation}
\label{forms}
p^*\eta=\sum_{a,b} f_{ab}\omega^{\a_1}\wedge\ldots\wedge
\omega^{\a_r}\wedge
\ov{\omega}^{\b_1}\wedge\ldots\wedge\ov{\omega}^{\b_s}
\end{equation}
on $K$, with coefficients 
$f_{ab}\hin C^{\infty}(K)$. 
Conversely, every $V$-invariant element of $C^{\infty}(K)\otimes
\bigwedge\tau(\n)^*\otimes\bigwedge\n^*$ is the pullback to $K$
of a differential form on $X$.
A form $\eta$ on $X$
is of $(p,q)$ type precisely when every summand on the right hand side of
(\ref{forms}) involves $p$ unbarred and $q$ barred terms. 

\begin{defn}
\label{invdefn}
{\em Inv}${}_{\R}(X)$ (respectively {\em Inv}${}_{\Q}(X)$)
denotes the ring of $K$-invariant forms in the 
$\R$-subalgebra (respectively $\Q$-subalgebra) of $A(X)$ generated by 
$
\{\frac{i}{2\pi}\om^{\a}\wedge\ov{\om}^{\b} \ | \ \a,\b \hin \Psi \}.
$
\end{defn}

Suppose now that $\pi:V\ra GL(E_0)$ is an irreducible unitary representation
of $V$ on a complex vector space $E_0$. $\pi$ defines a homogeneous
vector bundle $\ov{E}=K\times_VE_0\ra X$ which has a $K$-invariant hermitian 
metric. Extend $\pi$ to a unique holomorphic representation
$\pi:P\ra GL(E_0)$, and denote the induced representation of $\p$ by the
same letter. Then $\ov{E}=G\times_PE_0$ is a holomorphic hermitian 
vector bundle over $X$ which gives a complex structure to $K\times_VE_0$.

In [GrS], equation $(4.4)_X$, Griffiths and Schmid 
calculate the $K$-invariant curvature matrix 
$K(\ov{E})$ explicitly in terms of the above data. Their result is
\begin{equation}
\label{grs}
 K(\ov{E})=\sum_{\a,\b\in\Psi}\pi([e_{\a},e_{-\b}]_{\v})\otimes
\omega^{\a}\wedge\ov{\omega}^{\b}. 
\end{equation}
The invariant differential forms giving the Chern classes of homogeneous
line bundles were given by Borel [B]; see the introduction to [GrS]
for more references.

  Let $Y=F(\C)\simeq SL(n,\C)/B=SU(n)/S(U(1)^n)$ be the complex flag variety
and let $\ov{E}$ denote
the trivial hermitian vector bundle over $Y$, with the
tautological hermitian filtration
\[
\dis
\ov{\E}:\ 0 \subset \ov{E}_1 \subset \ov{E}_2 \subset\cdots 
\subset \ov{E}_n=\ov{E}
\]
with quotient line bundles $\ov{L}_i$ and all metrics induced from the
metric on $\ov{E}$. All of these bundles are homogeneous, and we want
to use equation (\ref{grs}) to compute their curvature matrices.
Note that (\ref{grs}) applies directly only to the line bundles
$\ov{L}_i$, as the higher rank bundles are not given by irreducible 
representations of the torus $S(U(1)^n)$. We can avoid this problem
by considering the grassmannian $Y_k=G_k(\C)=SU(n)/S(U(k)\times U(n-k))$
 and the natural projection $\rho:Y\ra Y_k$. Now (\ref{grs}) can be applied 
to the universal bundle $\ov{E}_k$ over $Y_k$
and the curvature matrix $K(\ov{E}_k)$ pulls back via $\rho$ to
the required matrix over $Y$. In fact by projecting to a
partial flag variety one can compute the curvature matrix of any
quotient bundle $E_l/E_k$. 
The representations $\pi$ of $V$ inducing these bundles are the
obvious ones in each case. What remains is a straightforward application of
equation (\ref{grs}),
 so we will describe the answer without belaboring the details.

 We have defined differential forms $\om^{ij}$, $\ov{\om}^{ij}$ on
$K=SU(n)$ which we identify with corresponding forms on $Y$. With
this notation, we can state (compare [GrS], $(4.13)_X$) :
\begin{prop}
\label{curvmat}
Let $k<l$ and consider the vector bundle $Q_{lk}=E_l/E_k$ over $F(\C)$. Let
the curvature matrix of $\ov{Q}_{lk}$ with its induced metric be
$\Theta=\{\Theta_{\a\b}\}_{k+1\lequ \a,\b\lequ l}$. Then
\[
\dis
\Theta_{\a\b}=\sum_{i\lequ k}\omega^{\a i}\wedge\ov{\omega}^{\b i}-
\sum_{j>l}\omega^{j\a}\wedge\ov{\omega}^{j\b}.
\]
\end{prop}

For notational convenience we let $\om_{ij}:=\gamma\om^{ji}$,
$\ov{\om}_{ij}:=\gamma\ov{\om}^{ji}$ and 
$\Om_{ij}:=\om_{ij}\wedge\ov{\om}_{ij}$,\ $(i<j)$,
where $\gamma$ is a constant such that 
$\gamma^2=\frac{i}{2\pi}$. Then we have
\begin{cor}
\label{grscor1}
\[
\dis
c_1(\ov{L}_k)=\sum_{i<k}\Om_{ik}-\sum_{j>k}\Om_{kj}
\]
\[
\dis
K_{\ov{E}_k}=\frac{i}{2\pi}K(\ov{E}_k)=-\left\{\sum_{j>k}\om_{\a j}
\wedge\ov{\om}_{\b j}\right\}_{1\lequ \a,\b \lequ k}
\]
\end{cor}

\medskip

Let $\dis \Omega:=\bigwedge_{i<j}\Om_{ij}$.
To compute classical
intersection numbers on the flag variety using the differential
forms in Corollary \ref{grscor1}
it suffices to know $\dis\int_Y\Omega$. If $x_i=-c_1(\ov{L}_i)$, it
is well known that $\eta=\S_{w_0}(x)=x_1^{n-1}x_2^{n-2}\cdots x_{n-1}$
is dual to the
class of a point in $Y$; thus $\dis\int_Y\eta=1$. An easy calculation shows
that $\dis\eta=\prod_{k=1}^{n-1}k!\cdot \Omega$, thus $\dis
\int_Y\Omega=\prod_{k=1}^{n-1}\frac{1}{k!}$.

\section{Invariant arithmetic Chow rings}
\label{afvs}

 It is well known that the arithmetic variety
$F$ has a  cellular decomposition in the sense 
of [Fu1], Ex. 1.9.1. It follows that one can use the excision exact sequence
for the groups $CH^{*,*}(F)$ 
(cf. [G], \sec 8) to show that $CH^{p,p-1}(F)=0$ (compare [Ma], Lemma 4.0.6).
 Therefore the exact sequence (\ref{ex1}) summed over $p$ gives
\begin{equation}
\label{flagex}
 0 \longrightarrow \wt{A}(F_{\R})
\stackrel{a}\longrightarrow \wh{CH}(F) 
\stackrel{\zeta}\longrightarrow CH(F)\longrightarrow 0.
\end{equation}
Recall that $\wt{A}(F_{\R})=\Ker\zeta$ is an ideal of $\wh{CH}(F)$ whose
$\wh{CH}(F)$-module structure is given as follows: if $\a\hin\wh{CH}(F)$ and
$\eta\hin\wt{A}(F_{\R})$, then $\a\cdot\eta=\om(a)\wedge\eta$. $\wt{A}(F_{\R})$
is not a square zero ideal, but its product is induced by $\theta\cdot
\eta=(dd^c\theta)\wedge\eta$. This product is well defined and commutative
([GS1], 4.2.11).

We equip $E(\C)$ with a trivial hermitian metric.
This metric induces metrics on all
the $L_i$, which thus become hermitian line bundles $\ov{L}_i$. 
Recall from \sec \ref{classgps} that $CH(F)$ 
has a free $\Z$-basis of monomials in the Chern classes $c_1(L_i)$. 
The unique map of abelian groups $\epsilon:CH(F)\ra \wh{CH}(F)$ 
sending $\prod c_1(L_i)^{k_i}$ to $\prod \wh{c}_1(\ov{L}_i)^{k_i}$
when $k_i\lequ n-i$ for all $i$ is then
a splitting of (\ref{flagex}). Thus we have an isomorphism of abelian groups
\begin{equation}
\label{bigiso}
\wh{CH}(F)\simeq CH(F)\oplus \wt{A}(F_{\R}).
\end{equation}

\medskip

As an analogue of the Arakelov Chow ring we define an {\em invariant
arithmetic Chow ring} $\wh{CH}_{inv}(F)$ as follows. 
Let $\Inv^{p,p}(F_{\R})$ be the group of $(p,p)$-forms $\eta$ 
in $\Inv_{\R}(F(\C))$ satisfying $F_{\infty}^*\eta=(-1)^p\eta$, and
set $\Inv(F_{\R})=\oplus_p \Inv^{p,p}(F_{\R})$. 
Let $\wt{\Inv}(F_{\R})\subset \wt{A}(F_{\R})$ be the image of $\Inv(F_{\R})$
in $\wt{A}(F_{\R})$.  Define the rings $\Inv(F_{\Q})$ and 
$\wt{\Inv}(F_{\Q})$ similarly, replacing $\R$ by $\Q$ in the above.

\begin{defn}
\label{invdef}
The invariant arithmetic Chow ring 
$\wh{CH}_{inv}(F)$ is the subring of $\wh{CH}(F)$ generated by
$\epsilon(CH(F))$ and $a(\wt{\mbox{\em Inv}}(F_{\R}))$.
\end{defn}

Suppose that $x,y\hin CH(F)$ and view $x$ and $y$ as elements of $\wh{CH}(F)$
using the inclusion $\epsilon$. In \sec \ref{ai}
 we will see that under the isomorphism (\ref{bigiso}),
$xy\hin CH(F)\oplus \wt{\Inv}(F_{\Q})$.
 It follows that there is an exact sequence of abelian groups
\begin{equation}
\label{invex}
0 \longrightarrow \wt{\Inv}(F_{\R}) 
\stackrel{a}\longrightarrow \wh{CH}_{inv}(F) 
\stackrel{\zeta}\longrightarrow CH(F)\longrightarrow 0
\end{equation}
which splits as before, giving
\begin{thm}
\label{chinv}
There is an isomorphism of abelian groups 
\[
\wh{CH}_{inv}(F)\simeq CH(F)\oplus\wt{\mbox{\em Inv}}(F_{\R}).
\]
\end{thm}

\medskip

\noindent
{\bf Remark 1:}
 One can define another `invariant arithmetic Chow ring' 
\[
\wh{CH}_{inv}^{\pr}(F):=\om^{-1}(\Inv(F_{\R})),
\]
where $\om$ is the ring homomorphism defined in \sec \ref{ait}.
There is a 
natural inclusion $\wh{CH}_{inv}(F)\hookrightarrow\wh{CH}_{inv}^{\pr}(F)$;
we do not know if these two rings coincide.

\medskip

\noindent
{\bf Remark 2:}
The arithmetic Chern classes of the
natural homogeneous vector bundles over $F$ are all contained in the ring
 $\wh{CH}_{inv}(F)$. In fact one need not use real coefficients for this;
it suffices to take $CH(F)\oplus \wt{\Inv}(F_{\Q})$ with the induced product
from $\wh{CH}(F)$. As there are bounds on the denominators that occur,
it follows 
 that {\em the subring of $\wh{CH}(F)$ generated by $\epsilon(CH(F))$
is a finitely generated abelian group}. However it seems that this group is
too small to contain the characteristic classes of all the vector bundles
of interest.

\section{Calculating arithmetic intersections}
\label{ai}

 In this section we describe an effective procedure for computing
arithmetic intersection numbers on the complete flag variety $F$.
One has a tautological
hermitian filtration of the trivial bundle $\ov{E}$ over $F$
\[
\dis
\ov{\E}:\ 0 \subset \ov{E}_1 \subset \ov{E}_2 \subset\cdots 
\subset \ov{E}_n=\ov{E}
\]
as in \sec \ref{hvb}. 
Recall 
that the inverse of the isomorphism
$CH(F)\simeq P_n/I_n$ sends $[X_i]$ to $-c_1(L_i)$. 
Let $x_i=-c_1(\ov{L}_i)$ and
 $\wh{x}_i=-\wh{c}_1(\ov{L}_i)$ for $1\lequ i \lequ n$.

If $\phi\hin \Lambda_n\otimes_{\Z}\Q$ is a homogeneous 
symmetric polynomial of positive degree then $\phi$ defines
a characteristic class. Theorem \ref{abc} applied to the
hermitian filtration $\ov{\E}$ shows that
\[
\dis
\phi(\wh{x}_1,\wh{x}_2,\ldots,\wh{x}_n)=(-1)^{\deg \phi}\,\wt{\phi}(\ov{\E})
\]
in the arithmetic Chow ring $\wh{CH}(F)$.
In particular for $\phi=e_i$ an elementary symmetric polynomial this gives
\[
\dis
e_i(\wh{x}_1,\wh{x}_2,\ldots,\wh{x}_n)=(-1)^i\,\wt{c}_i(\ov{\E}).
\]

Let $h$ be a homogeneous polynomial in the ideal $I_n$. We will give an
algorithm for computing $h(\wh{x}_1,\wh{x}_2,\ldots\wh{x}_n)$ as a class 
in $\wt{\Inv}(F_{\Q})$:

\medskip

\noindent
{\bf\underline{Step 1}:} Decompose $h$ as a sum $h=\sum e_if_i$ for some polynomials 
$f_i$. More canonically one may use the equality
\[
\dis
h=\sum_{w\in S_n} \left< h , \S^w \right> \S_w
\]
from \sec \ref{classgps}. 
Since $a(x)y=a(x\omega(y))$ in $\wh{CH}(F)$
and $\omega(f_i(\wh{x}_1,\ldots,\wh{x}_n))=
f_i(x_1,\ldots,x_n)$, we have
\[
\dis
h(\wh{x}_1,\wh{x}_2,\ldots\wh{x}_n)=\sum_{i=1}^n
(-1)^i\,\wt{c}_i(\ov{\E})f_i(x_1,x_2,\ldots,x_n)=
\]
\[
\dis
=\sum_{w\in S_n} (-1)^{\deg h +l(w)}\,\wt{\left< h , \S^w \right>}(\ov{\E}) 
\S_w(x_1,x_2,\ldots,x_n).
\]

\noindent
{\bf\underline{Step 2}:}
By Corollary \ref{ratbcf}, we may express 
the forms $\wt{c}_i(\ov{\E})$ and $\dis \wt{\left< h , \S^w \right>}(\ov{\E})$
as polynomials
in the entries of the matrices $K_{E_i}$ and $K_{L_i}=c_1(\ov{L_i})$ with
rational coefficients. 
In practice this may be done recursively for the Chern forms $\wt{c}_i$
as follows: Use equation (\ref{*}) and
the construction in Corollary \ref{ratcor} to
obtain the power sum forms $\wt{p}_i(\ov{\E})$, then apply the formulas
(\ref{sumprod}) to Newton's identity
\[
\dis
p_i-c_1p_{i-1}+c_2p_{i-2}-\cdots+(-1)^iic_i=0.
\]
On the other hand 
Corollary \ref{grscor1} gives explicit expressions for all
the above curvature matrices in terms of differential forms on $F(\C)$.
 Thus we obtain formulas for 
$\wt{c}_i(\ov{\E})$ and $\dis \wt{\left< h , \S^w \right>}(\ov{\E})$
in terms of these forms.
For example, using the notation of \sec \ref{hvb}, we have
\begin{prop}
$\wt{c}_1(\ov{\E})=0$ and $\dis\wt{c}_2(\ov{\E})=-\sum_{i<j}\Om_{ij}$.
\end{prop}
{\bf Proof.} Use (\ref{*}), properties (a) and (b) at the end of 
\sec \ref{cbcfs}, and the identity $2c_2=c_1^2-p_2$.\ \ \
\endproof

\medskip

\noindent
{\bf\underline{Step 3}:}
Substitute the forms obtained in Step 2 into the formulas 
given in Step 1. Note that
the result is the class of a form in $\Inv(F_{\Q})$ since all the
ingredients are functorial for the natural $U(n)$ action on $F(\C)$.

\medskip

 In particular, if $k_i$ are nonnegative integers with $\sum k_i=
\dim{F}={n \choose 2}+1$, the monomial $X_1^{k_1}\cdots X_n^{k_n}$
is in the ideal $I_n$. 
If $X_1^{k_1}\cdots X_n^{k_n}=\sum e_if_i$, 
 then we have 
\[
\dis
\wh{x}_1^{k_1}\wh{x}_2^{k_2}\cdots\wh{x}_n^{k_n}=
\sum_i(-1)^i\,\wt{c}_i(\ov{\E})f_i(x_1,\ldots,x_n).
\]
Now if $\Om=\bigwedge \Om_{ij}$ is defined as in \sec \ref{hvb}, we have
shown that
\[
\dis
\wt{c}_i(\ov{\E})f_i(x_1,\ldots,x_n)=r_i\Om
\]
for some rational number $r_i$. Therefore 
\[
\dis
\wh{\deg}(\wh{x}_1^{k_1}\wh{x}_2^{k_2}\cdots\wh{x}_n^{k_n})=
\frac{1}{2}\sum_i(-1)^i\,r_i\int_{F(\C)}\Om=
\frac{1}{2}\sum_i(-1)^i\,r_i\prod_{k=1}^{n-1}\frac{1}{k!}\, .
\]
Of course this equation implies
\begin{thm}
\label{mainthm} The arithmetic Chern number
$
\dis
\wh{\deg}(\wh{x}_1^{k_1}\wh{x}_2^{k_2}\cdots\wh{x}_n^{k_n})
$
is a rational number.
\end{thm}

\noindent
{\bf Remark.}\ 
For $a<b$ let $\ov{Q}_{b,a}=E_b/E_a$, equipped with the induced metric.
 Then one can show that any intersection number
$
\dis
\wh{\deg}(\prod_i \wh{c}_{m_i}(\ov{Q}_{b_ia_i})^{k_i})$ for
$
\sum k_im_i(b_i-a_i)=\dim F
$
is rational. This is done by using 
the hermitian filtrations
\[
\dis
0\subset \ov{Q}_{a+1,a}\subset\ov{Q}_{a+2,a}\subset\cdots\subset\ov{Q}_{b,a}
\]
and Theorem \ref{abc} to reduce the problem to the intersections occuring in 
Theorem \ref{mainthm}. To compute arithmetic intersections of the form
$(0,\eta)\cdot (0,\eta^{\pr})$ with $\eta,\eta^{\pr}\hin \wt{\Inv}(F_{\Q})$,
we need to know the value of $dd^c\eta$. For this one may use the 
Maurer-Cartan structure equations on $SU(n)$ (cf. [GrS], Chp. 1); all
such intersections lie in $\wt{\Inv}(F_{\Q})$.

\medskip
 
 Although there is an effective algorithm for computing arithmetic 
Chern numbers,
explicit general formulas seem difficult to obtain.
 There are some
general facts we can deduce for those intersections that pull back from
grassmannians, for instance that $\wh{x}_1^{n+1}=\wh{x}_n^{n+1}=0$.
There is also a useful symmetry in these intersections:
\begin{prop}
\label{symprop}
$
\wh{x}_1^{k_1}\wh{x}_2^{k_2}\cdots\wh{x}_n^{k_n}=
\wh{x}_n^{k_1}\wh{x}_{n-1}^{k_2}\cdots\wh{x}_1^{k_n},
$ \ for all integers $k_i\gequ 0$.
\end{prop}
{\bf Proof.} This is a consequence of the involution $\nu : F(\C)
\ra F(\C)$ sending 
\[
\dis
\ov{\E} : \ 0 \subset \ov{E}_1 \subset \ov{E}_2 \subset \cdots \subset
\ov{E}_n=\ov{E}
\]
to
\[
\dis
\ov{\E}^{\bot} : \ 0=\ov{E}^{\bot} \subset \ov{E}_{n-1}^{\bot} 
\subset \ov{E}_{n-2}^{\bot}
 \subset \cdots \subset 0^{\bot}=\ov{E}.
\]
Over $\Spec\Z$, $\nu$ corresponds to the map of flag varieties sending
$E_i$ to the quotient $E/E_i$.
If $\wh{x}_i^{\bot}$ are the arithmetic Chern classes obtained from
$\ov{\E}^{\bot}$, then using the split exact sequences
$0\ra \ov{E}_i \ra \ov{E} \ra \ov{E}_i^{\bot} \ra 0$ we obtain
\[
\dis
\wh{x}_i^{\bot}=-\wh{c}_1(\ov{L}_i^{\bot})=
-\wh{c}_1(\ov{E}_{n-i}^{\bot})+\wh{c}_1(\ov{E}_{n+1-i}^{\bot})=
\wh{c}_1(\ov{E}_{n-i})-\wh{c}_1(\ov{E}_{n+1-i})=\wh{x}_{n+1-i}.
\]
Since $\nu$ is an isomorphism, the result follows.
\endproof

\section{Arithmetic Schubert calculus}
\label{asc}

 Let $P_n$, $I_n$, $\Lambda_n$ and $S^{(n)}$ be as in \sec \ref{classgps}. 
The Chow ring $CH(F)$
is isomorphic to the quotient $H_n=P_n/I_n$. Recall that $H_n$ has a
natural basis of Schubert polynomials $\{\S_w\ | \ w\hin S_n \}$, 
and that the $\S_w$ for $w\hin S^{(n)}$ form a free $\Z$-basis of $P_n$.
We let $T_n=S^{(n)}\smallsetminus S_n$.
The key property of Schubert polynomials that we require for
the `arithmetic Schubert calculus' is described in 

\begin{lemma}
\label{schlemma}
If $w\hin T_n$, then $\S_w\hin I_n$. In fact we
have a decomposition
\[
\dis
\S_w=\sum_{v\in S_n}\left<\S_w,\S^v\right>\S_v,
\]
 where $ \left<\S_w,\S^v\right>\hin \Lambda_n\cap I_n$.
\end{lemma}
{\bf Proof.} Assume first that $w(1)>w(2)>\cdots >w(n)$, so that $w$
is {\em dominant}. Then by [M], (4.7) we have
\[
\dis
\S_w=X_1^{w(1)-1}X_2^{w(2)-1}\cdots X_n^{w(n)-1}.
\]
If $w\notin S_n$ then clearly $w(1)>n$, so $X_1^{w(1)-1}\hin I_n$ and
thus $\S_w\hin I_n$.

If $w\hin T_n$ is
arbitrary, form $w^{\pr}\hin T_n$ by rearranging
$(w(1),w(2),\ldots,w(n))$ in decreasing order and letting $w^{\pr}(i)=
w(i)$ for $i>n$. We have shown that $\S_{w^{\pr}}\hin I_n$. There is an
element $v\hin S_n$ such that $wv=w^{\pr}$ and $l(v)=l(w^{\pr})-l(w)$.
Note that since $\partial_v$ is $\Lambda_n$-linear,
$\partial_v I_n\subset I_n$.
Therefore  ([M], (4.2)):
$\S_w=\partial_v\S_{wv}=\partial_v\S_{w^{\pr}}\hin I_n.$

The decomposition claimed now follows, as in \sec \ref{classgps}.
\endproof

\medskip

It is well known that there is an equality in $P_{\infty}$
\begin{equation}
\label{cuvws}
\S_u\S_v=\sum_{w\in S_{\infty}}c_{uv}^w\S_w,
\end{equation}
where the $c_{uv}^w$ are nonegative integers that vanish whenever
$l(w)\neq l(u)+l(v)$\ ([M], (A.6)). 
A particular
case of this is {\em Monk's formula:} if $s_k$ denotes the transposition
$(k,k+1)$, then
\[
\dis
\S_{s_k}\S_w=\sum_t \S_{wt}
\]
summed over all transpositions $t=(i,j)$ such that $i\lequ k <j$
and $l(wt)=l(w)+1$\ ([M], ($4.15^{\pr\pr}$)).

\medskip

 We now express arithmetic intersections in $\wh{CH}(F)$ using the
basis of Schubert polynomials. Lemma \ref{schlemma} is
the main reason why this basis facilitates our task.
 This property (for Schur functions) also plays a 
crucial role in the arithmetic Schubert calculus for grassmannians
(see \sec \ref{pfvs} and [Ma], Th. 5.2.1).

  For each $w\hin S^{(n)}$, let $\wh{\S}_w=\S_w(\wh{x}_1,
\ldots,\wh{x}_n)$. If $w\hin T_n$
 then Lemma \ref{schlemma} and the discussion in \sec \ref{ai}
imply that $\wh{\S}_w\hin \wt{\Inv}(F_{\Q})$; we
denote these classes by $\wt{\S}_w$. We have
\[
\dis
\wt{\S}_w=\sum_{v\in S_n}(-1)^{l(v)+l(w)}
\wt{\left<\S_w,\S^v\right>}(\ov{\E})\S_v(x_1,\ldots,x_n).
\]
We can now describe the multiplication in $\wh{CH}_{inv}(F)$:
\begin{thm}
\label{slring}
Any element of $\wh{CH}_{inv}(F)$ can be expressed uniquely in the form
$\dis \sum_{w\in S_n}a_w\wh{\S}_w+\eta$, where $a_w\hin\Z$ and 
$\eta\in\wt{\mbox{\em Inv}}(F_{\R})$. For $u,v\hin S_n$ we have
\begin{equation}
\label{punchline}
\wh{\S}_u\cdot\wh{\S}_v=\sum_{w\in S_n}c_{uv}^w\wh{\S}_w+
\sum_{w\in T_n}c_{uv}^w\wt{\S}_w,
\end{equation}
\[
\dis
\wh{\S}_u\cdot \eta=\S_u(x_1,\ldots,x_n)\wedge\eta,
\ \ \ \ \mbox{and}
\ \ \ \ \eta\cdot \eta^{\pr}=(dd^c\eta)\wedge\eta^{\pr},
\]
where $\wt{\S}_w\in\wt{\mbox{\em Inv}}(F_{\Q})$, $\eta$,
$\eta^{\pr}\in\wt{\mbox{\em Inv}}(F_{\R})$ and the $c_{uv}^w$ are as in 
{\em (\ref{cuvws})}.
\end{thm}
{\bf Proof.} The first statement is a corollary of Theorem \ref{chinv}.
Equation (\ref{punchline}) follows immediately from the formal identity
(\ref{cuvws}) and our definition of $\wt{\S}_w$. The rest is a consequence
of properties of the multiplication in $\wh{CH}(F)$ discussed in 
\sec \ref{afvs} and \sec \ref{ai}. \endproof

\medskip

\noindent
{\bf Remark.}
It is interesting to note that we also have, for $u,v\hin T_n$, 
\[
\dis
\wt{\S}_u\cdot\wt{\S}_v=(dd^c\wt{\S}_u)\wedge\wt{\S}_v=
\sum_{w\in T_n}c_{uv}^w\wt{\S}_w
\]
in $\wt{\Inv} (F_{\Q})$. 

\medskip

 Applying (\ref{punchline}) when $\S_u=\S_{s_k}$ is a special Schubert class
gives

\begin{cor} (Arithmetic Monk Formula): 
\[
\dis
\wh{\S}_{s_k}\cdot\wh{\S}_w=\sum_s \wh{\S}_{ws} +
\sum_t \wt{\S}_{wt},
\]
where the first sum is over all transpositions $s=(i, j)\hin S_n$ such that 
$i\lequ k <j$ and $l(ws)=l(w)+1$, and the second over all transpositions
$t=(i,n+1)$ with $i\lequ k$ and 
 $w(i)>w(j)$ for all $j$ with $i<j \lequ n$.
\end{cor}

\section{Examples}
\label{ex}
\subsection{Heights}

The flag variety $F$ has a natural pluri-Pl\"{u}cker embedding
$j:F\hookrightarrow \P^N_{\Z}$.
$j$ is defined as the composition of a product of Pl\"{u}cker embeddings
followed by a Segre embedding; if $Q_i=E/E_i$,
then $j$ is associated to the
line bundle $\dis Q=\bigotimes_{i=1}^{n-1} \det (Q_i)$. 
Let $\ov{\O}(1)$ denote the
canonical line bundle over projective space, equipped with its canonical
metric (so that $c_1(\ov{\O}(1))$ is the Fubini-Study form). The 
{\em height} of $F$ relative to $\ov{\O}(1)$ (cf.
[Fa], [BoGS], [S]) is defined by 
\[
\dis
h_{\ov{\O}(1)}(F)
=\wh{\deg}(\wh{c}_1(\ov{\O}(1))^{{n \choose 2}+1}\vert \ F).
\]
Since
\[
\dis
j^*(\wh{c}_1(\ov{\O}(1)))=\wh{c}_1(\ov{Q})=
-\sum_{i=1}^{n-1}\wh{c}_1(\ov{E}_i)=
\sum_{i=1}^{n-1} (n-i)\wh{x}_i=
\sum_{i=1}^{n-1} \wh{\S}_{s_i},
\]
we have that 
\[
\dis
h_{\ov{\O}(1)}(F)=
\wh{\deg}(\wh{c}_1(\ov{Q})^{{n \choose 2}+1}\vert \ F)=
\wh{\deg}((\sum_{i=1}^{n-1} \wh{\S}_{s_i})^{{n \choose 2}+1}).
\]
Now Theorems \ref{mainthm} and \ref{slring} immediately imply
\begin{thm}
\label{slheight}
The height $h_{\ov{\O}(1)}(F)$ is a rational number. 
\end{thm}


\subsection{Intersections in $F_{1,2,3}$}

In this section we calculate the arithmetic intersection numbers for
the classes $\wh{x}_i$ in $\wh{CH}(F)$ when $n=3$, so $F=F_{1,2,3}$.

Over $F$ we have 3 exact sequences
\[
\dis
\ov{\E}_i : \ 0\ra \ov{E}_{i-1} \ra \ov{E}_i \ra \ov{L}_i \ra 0
\ \ \ \ \ \
1\lequ i \lequ 3.
\]
We adopt the notation of \sec \ref{hvb} and define $\Om_{ij}=
\om_{ij}\wedge\ov{\om}_{ij}$. Then Corollary \ref{grscor1} gives
\[
\dis
x_1=\Om_{12}+\Om_{13},
\ \ \ \ 
x_2=-\Om_{12}+\Om_{23},
\ \ \ \
x_3=-\Om_{13}-\Om_{23},
\]
\[
\dis
K_{E_2}=-\left(
\begin{array}{cc}
\Om_{13} & \om_{13}\wedge\ov{\om}_{23} \\
\om_{23}\wedge\ov{\om}_{13} & \Om_{23}
\end{array}\right).
\]
We refer now to the properties of the forms $\wt{c_k}$ mentioned at the
end of \sec \ref{cbcfs}.
By property (a) $\wt{c}(\ov{\E}_1)=0$, while (b) gives
$\wt{c}(\ov{\E}_2)=-\Om_{12}$. Property (c) applied to
 $\ov{\E}_3$ gives $\wt{c}(\ov{\E}_3)=-\Om_{13}-\Om_{23}+
3\Om_{13}\Om_{23}$.
Using the construction of the Bott-Chern form for the total Chern class
given in the proof of Theorem \ref{bcf}, we find that
\begin{equation}
\label{firstkey}
\wt{c}(\ov{\E})=
-\Om_{12}-\Om_{13}-\Om_{23}-\Om_{12}\Om_{13}-\Om_{12}\Om_{23}+
3\Om_{13}\Om_{23}.
\end{equation}
Notice that this expression for $\wt{c}(\ov{\E})$ is not unique as a 
class in $\wt{\Inv}(F_{\R})$.
For instance, we can add the exact form 
$c_1(\ov{L}_1)c_1(\ov{L}_2)-c_2(\ov{E}_2)=\Om_{12}\Om_{23}-
\Om_{12}\Om_{13}-\Om_{13}\Om_{23}$ to get
\begin{equation}
\label{key}
\wt{c}(\ov{\E})=
-\Om_{12}-\Om_{13}-\Om_{23}-2\Om_{12}\Om_{13}+2\Om_{13}\Om_{23}.
\end{equation}
The Bott-Chern form (\ref{key}) is the key to computing
 any intersection number $\wh{\deg}(\wh{x}_1^{k_1}
\wh{x}_2^{k_2}\wh{x}_3^{k_3})$, following the prescription of \sec
\ref{ai}. (Of course we can just as
well use (\ref{firstkey}), with the same results.)
 For example, since $x_1^4=x_1^3e_1-x_1^2e_2+x_1e_3$, we have
\[
\dis
\wh{x}_1^4=x_1^2(\Om_{12}+\Om_{13}+\Om_{23})+
x_1(2\Om_{12}\Om_{13}-2\Om_{13}\Om_{23})=
2\Om-2\Om=0.
\]
On the other hand, a similar calculation for $\wh{x}_2^4$ gives
\[
\dis
\wh{x}_2^4=-x_2^2\wt{c}_2(\ov{\E})-
x_2\wt{c}_3(\ov{\E})=-2\Om+4\Om=2\Om.
\]
Thus $\dis\wh{\deg}(\wh{x}_2^4)=\int_{F(\C)}\Om=\frac{1}{2}$.

The following is a table of all the intersection numbers  
$\wh{\deg}(\wh{x}_1^{k_1}\wh{x}_2^{k_2}\wh{x}_3^{k_3})$
(multiplied by 4):

\begin{center} \begin{tabular}{|cc|cc|cc|} \hline
&&&&& \\
$k_1k_2k_3$ & $4\,\wh{\deg}$ & $k_1k_2k_3$ & $4\,\wh{\deg}$ &
$k_1k_2k_3$ & $4\,\wh{\deg}$ \\ \hline
400 & 0 & 004 & 0 & 040 & 2 \\
310 & 5 & 013 & 5 & 121 & 2 \\
301 & -5 & 103 & -5 & 202 & 9 \\
220 & -1 & 022 & -1 && \\
211 & -4 & 112 & -4 && \\
130 & -1 & 031 & -1 && \\ \hline
\end{tabular} \end{center}
Note that the numbers in the first two columns are equal, in agreement with 
Proposition \ref{symprop}.
We can use the table to compute the height of $F$ in its pluri-Pl\"{u}cker
embedding in $\P^8_{\Z}$:
\[
\dis
h_{\ov{\O}(1)}(F_{1,2,3})=\wh{\deg}((2\wh{x}_1+\wh{x}_2)^4)=\frac{65}{2}.
\]

\section{Partial flag varieties}
\label{pfvs}

In this final section we show how to generalize the previous work
to partial flags $F(\r)$. Our results may thus be regarded as an
extension of those of Maillot [Ma] in the grassmannian case.

 As usual we have a tautological filtration of type $\r$
\begin{equation}
\label{fil3}
\E:\ 0 \subset E_1 \subset E_2 \subset\cdots \subset E_m=E
\end{equation}
of the trivial bundle over $F(\r)$, with quotient bundles $Q_i$.
 Equip $E(\C)$ with the trivial hermitian metric, inducing metrics on all the
above bundles. 
 The calculations of \sec \ref{hvb} apply equally well to $X_{\r}=F(\r)(\C)$. 
Proposition \ref{curvmat} describes the 
curvature matrices of all the relevant homogeneous vector bundles, and 
one can compute classical intersection 
numbers on $X_{\r}$ in a similar fashion.

 Call a permutation $w\hin S_{\infty}$ an
{\em $\r$-permutation} if $w(i)<w(i+1)$ for all $i$ not contained in 
$\{r_1,\ldots,r_{m-1}\}$. Let $S_{n,\r}$ and $T_{n,\r}$ be the set of
$\r$-permutations in $S_n$ and $T_n$, respectively.
For such $w$  one knows (cf. [Fu2], \sec 8)
that the Schubert polynomial
$\S_w$ is symmetric in the variables in each of the groups
\begin{equation}
\label{vars}
X_1,\ldots,X_{r_1};X_{r_1+1},\ldots,X_{r_2};\ldots ;
X_{r_{m-2}+1},\ldots,X_{r_{m-1}}.
\end{equation}

The product group $\dis H=\prod_{i=1}^m S_{r_i-r_{i-1}}$ 
acts on $P_n$, the factors for $i<m$
by permuting the variables in the corresponding group of (\ref{vars}), while
$S_{n-r_{m-1}}$ permutes the remaining variables $X_{r_{m-1}+1},\ldots ,X_n$.
If $P_{n,\r}=P_n^H$ is the ring of invariants and $I_{n,\r}=P_{n,\r}
\cap I_n$, then $CH(F(\r))\simeq P_{n,\r}/I_{n,\r}$.
The set of Schubert polynomials $\S_w$ for all $w\hin S_{n,\r}$
is a free $\Z$-basis for $P_{n,\r}/I_{n,\r}$.

Let $w\hin S_{n,\r}$.
If we regard each of the groups of variables (\ref{vars})
 as the Chern roots of the
bundles $Q_1,Q_2,\ldots,Q_{m-1}$, it follows that we may write $\S_w$ as a
polynomial $\S_{w,\r}$ in the Chern classes of the $Q_i$, 
$1\lequ i \lequ m-1$. The class of $\S_{w,\r}$ in $CH(F(\r))$
is that of corresponding Schubert variety in $F(\r)$ (see loc. cit. for
the relative case). By putting `hats' on all the quotient bundles involved
(with their induced metrics as in \sec \ref{asc}) we obtain classes 
$\wh{\S}_{w,\r}$ in $\wh{CH}_{inv}(F(\r))$.

 The analysis of \sec \ref{afvs} remains valid; the map $\epsilon$ can be
defined by $\epsilon(\S_{w,\r})=\wh{\S}_{w,\r}$.  In particular we 
have an invariant arithmetic Chow
ring $\wh{CH}_{inv}(F(\r))$ for which Theorem \ref{chinv} holds. If
$F(\r)=G_d$ is a Grassmannian over $\Spec\Z$, then $\wh{CH}_{inv}(G_d)$ 
coincides with the Arakelov Chow ring $CH(\ov{G_d})$, where $G_d(\C)$ is
given its natural invariant K\"{a}hler metric, as in [Ma].

 Suppose that $\r^{\pr}$ is a refinement of $\r$, so we have a projection
$p:F(\r^{\pr})\ra F(\r)$. In this case there are natural inclusions
$\wt{\Inv}(F(\r)_{\R})\hookrightarrow \wt{\Inv}(F(\r^{\pr})_{\R})$ and
$CH(F(\r))\hookrightarrow CH(F(\r^{\pr}))$.
Applying the five lemma to the two exact sequences (\ref{invex})
shows that the pullback
$p^*:\wh{CH}_{inv}(F(\r))\hookrightarrow \wh{CH}_{inv}(F(\r^{\pr}))$ is 
an injection. Note however that this is not compatible with the splitting
of Theorem \ref{chinv}.

 One can compute arithmetic intersections in $\wh{CH}_{inv}(F(\r))$ as in
\sec \ref{ai}. Applying Theorem \ref{abc} to the filtration (\ref{fil3})
(with induced metrics as above)
gives the key relation required for the calculation. In particular we see that
all the arithmetic Chern numbers are rational, as is the Faltings 
height of $F(\r)$ in its natural pluri-Pl\"{u}cker embedding. Theorem
\ref{slheight} thus generalizes the corresponding result of
Maillot mentioned in \sec \ref{intro}.

 There is an arithmetic Schubert calculus in $\wh{CH}_{inv}(F(\r))$ analogous
to that for complete flags. The analogue of Lemma \ref{schlemma} 
is true, that is $\S_w\hin I_{n,\r}$ if $w\hin T_{n,\r}$ 
(this is an easy consequence of Lemma \ref{schlemma} itself). 
It follows that for $w\hin T_{n,\r}$, $\wh{\S}_{w,\r}$ is a class
$\wt{\S}_{w,\r}\hin \wt{\Inv}(F(\r)_{\Q})$.
The analogue of (\ref{punchline}) in this context is
\begin{equation}
\label{partialpunch}
\wh{\S}_{u,\r}\cdot\wh{\S}_{v,\r}=\sum_{w\in S_{n,\r}}c_{uv}^w\wh{\S}_{w,\r}+
\sum_{w\in T_{n,\r}}c_{uv}^w\wt{\S}_{w,\r}
\end{equation}
where $u,v\hin S_{n,\r}$ and the numbers $c_{uv}^w$ are as in (\ref{cuvws}).
The remaining statements of Theorem \ref{slring} require no further change.

\medskip

\noindent
{\bf Remark.} Equation (\ref{partialpunch}) is not a direct generalization
of the analogous statement in [Ma], Theorem 5.2.1. However one can reformulate
Maillot's results using the classes $\wh{c}_*(\ov{S})$ instead of
$\wh{c}_*(\ov{Q}-\ov{\E})$ (notation as in [Ma], \sec 5.2). 
With this modification,
the arithmetic Schubert calculus described above (for $m=2$) and that in [Ma] 
coincide. In the grassmannian case $\S_{w,\r}$ is a Schur polynomial 
 and there are explicit formulas
for $\wt{\S}_{w,\r}$ in terms of harmonic forms on $G_d(\C)$ (as in loc. cit.).

\end{document}